%% file: Submitted.tex
\documentclass[11pt,a4paper]{article}
\input{Preamble.tex}

\ifPDFTeX\else
\fi
\newcommand{\F}{\mathcal{F}}
\newcommand{\spc}{\mathrm{sp}}
\newcommand{\wsd}{\mathrm{wsd}}
\newcommand{\DS}{\mathrm{DelSave}}
\newcommand{\Del}{\mathrm{Delete}}
\newcommand{\pa}{\partial}
\begin{document}

\title{Every IC-planar graph has weak-strict degeneracy at most $6$}
\author{
    Zhen He\thanks{School of Mathematics and Statistics, Henan University, Kaifeng, 475004, P.R. China.}
    \and
    Yilong Li\thanks{School of Mathematical Sciences, Zhejiang Normal University, Jinhua, P.R. China}
    \and
    Tao Wang\thanks{Center for Applied Mathematics, Henan University, Kaifeng, 475004, P.R. China. Corresponding author. Email: \texttt{wangtao@henu.edu.cn}; \href{https://orcid.org/0000-0001-9732-1617}{ORCID: 0000-0001-9732-1617}.}
    \and
    Xiaojing Yang\thanks{School of Mathematics and Statistics, Henan University, Kaifeng, 475004, P.R. China.}
}
\date{November 6, 2026}
\maketitle

\begin{abstract}
We prove that every IC-planar graph has weak-strict degeneracy at most $6$. Weak-strict degeneracy is the weak degeneracy of Bernshteyn and Lee plus one, so it bounds the coloring parameters controlled by weak degeneracy. In particular, every IC-planar graph is $6$-choosable, $6$-DP-colorable and $6$-paintable, which strengthens the theorem of Yang, Wang, Wang and Lih that IC-planar graphs are $6$-choosable. An IC-planar graph is a subgraph of a plane graph in which a set of pairwise vertex-disjoint inner $4$-faces, called marked faces, receive both diagonals. Our main tool is a Thomassen-type statement for such graphs: vertices of marked faces get one extra unit, except those in a prescribed independent set of boundary vertices. With no marked faces, this statement is exactly the theorem of Bernshteyn, Lee and Smith-Roberge that gives weak-strict degeneracy at most $5$ for planar graphs, so it contains both results. We prove it by a minimal-counterexample argument following Bernshteyn, Lee and Smith-Roberge. The new ingredients are a deletion rule that turns destroyed marked faces into ordinary edges, a pseudo-$2$-chord lemma for paths that use one diagonal, and a choice of orientation of the outer cycle that avoids a vertex of the independent set joined to two vertices of the prescribed boundary path.


\end{abstract}

\section{Introduction}
Bernshteyn and Lee~\cite{MR4606413} introduced the \emph{weak degeneracy} $\mathrm{wd}(G)$ of a graph $G$, a variant of degeneracy based on a modified greedy coloring procedure. We work with the following version, shifted by one. For a graph $G$ and $f:V(G)\to\mathbb Z$, the operation $\Del(G,f,u)$ removes $u$ and lowers $f$ by one on every neighbor of $u$; $\DS(G,f,u,w)$ (for a neighbor $w$ of $u$ with $f(u) > f(w)$) does the same except that $f(w)$ is not lowered. An operation is \emph{legal} if all resulting values are $\geq 1$. A graph $G$ is \emph{weak-strictly $f$-degenerate} if $f \geq 1$ and all vertices can be removed by legal operations, and \emph{weak-strictly $k$-degenerate} if this holds for $f \equiv k$. The \emph{weak-strict degeneracy} $\wsd(G)$ is the least $k$ such that $G$ is weak-strictly $k$-degenerate.

Both operations, and the condition $f(u) > f(w)$, commute with adding a constant to $f$. Hence $G$ is weak-strictly $f$-degenerate exactly when it is weakly $(f-1)$-degenerate in the sense of Bernshteyn and Lee~\cite{MR4606413}, where legality means values $\geq 0$, and $\wsd(G) = \mathrm{wd}(G) + 1$; the name weak-strict degeneracy follows~\cite{MR5077932}. If $G$ is weakly $d$-degenerate, then the chromatic number, the list chromatic number, the DP-chromatic number and the paint number of $G$ are all at most $d + 1$~\cite{MR4606413}. Thus $\wsd(G)$ is an upper bound for all of these coloring parameters. Throughout the paper we work with weak-strict degeneracy.

Thomassen~\cite{MR1290638} proved that every planar graph is $5$-choosable. His proof rests on a stronger statement about near-triangulations with a precolored boundary path and lists of size $3$ on the other boundary vertices. Bernshteyn, Lee and Smith-Roberge~\cite{MR5126524} proved the weak-degeneracy analog: planar graphs are weakly $4$-degenerate, that is, $\wsd(G) \leq 5$ for every planar graph $G$. This corrects an error in the planar case of~\cite{MR4606413}, and the main tool is a Thomassen-type statement \cite[Theorem~2.2]{MR5126524}. Weak degeneracy of planar graphs and related classes was also studied in~\cite{Davies2025,MR4663366,MR4564473}.

A graph is \emph{$1$-planar} if it can be drawn in the plane so that each edge is crossed at most once, and \emph{IC-planar} (independent crossings) if, in addition, no two crossed edge pairs share an end vertex. Albertson~\cite{MR2434266} introduced IC-planar graphs and conjectured that they are $5$-colorable; this was proved by Kr\'al' and Stacho~\cite{MR2674491}. Every IC-planar graph on $n$ vertices has at most $3.25n-6$ edges~\cite{MR3000647}, so IC-planar graphs are $6$-degenerate, and hence $\wsd(G) \leq 7$. Yang, Wang, Wang and Lih~\cite{MR4296752} proved that every IC-planar graph is $6$-choosable. Their proof also follows Thomassen: they replace each crossing by an inner $4$-face, so that the $4$-faces are pairwise vertex-disjoint, and prove a list-coloring statement in which the vertices of each inner face receive distinct colors \cite[Theorem~1.2]{MR4296752}. For list coloring of graphs with few or sparse crossings, see also~\cite{MR2873215,MR3597095}.

Our main result strengthens the theorem of Yang et al.\ to weak-strict degeneracy.

\begin{theorem}\label{thm:main}
Every IC-planar graph has weak-strict degeneracy at most $6$. Consequently, every IC-planar graph is $6$-choosable, $6$-DP-colorable and $6$-paintable.
\end{theorem}

We use the same reduction to plane graphs with pairwise vertex-disjoint marked $4$-faces as~\cite{MR4296752}, but the diagonals of the marked faces are now edges of the graph. The main tool, \zcref{thm:C} in \zcref{sec:inductive}, extends \cite[Theorem~2.2]{MR5126524} to this setting: with no marked faces it is exactly that theorem. A list-coloring statement such as \cite[Theorem~1.2]{MR4296752} does not imply a weak-degeneracy statement, so the proof of \zcref{thm:C} follows the minimal-counterexample argument of~\cite{MR5126524} rather than~\cite{MR4296752}.

\section{Preliminaries}
We need three basic properties of weak-strict degeneracy.

\begin{lemma}[{\cite[Lemma~2.1]{MR4606413}}]\label{lem:mono}
If $G$ is weak-strictly $f$-degenerate and $g \geq f$ pointwise, then $G$ is weak-strictly $g$-degenerate.
\end{lemma}
\begin{proof}
This is \cite[Lemma~2.1]{MR4606413} applied to $f-1 \leq g-1$.
\end{proof}

\begin{lemma}[subgraphs]\label{lem:sub}
If $G$ is weak-strictly $f$-degenerate and $G'\subseteq G$ is a subgraph, then $G'$ is weak-strictly $f|_{V(G')}$-degenerate.
\end{lemma}
\begin{proof}
It suffices to treat $G' = G - e$ and $G' = G - x$; induct on $|V(G)|$. Let the first operation on $G$ remove $u$ (saving $w$, or nobody), producing $(G - u, f_1)$, weak-strictly $f_1$-degenerate.

\emph{Edge $e$.} If $u \notin e$, perform the same operation on $G - e$: it is legal, and the result is $(G - u - e, f_1)$, done by induction. If $e = uy$, perform the same operation on $G - e$ (if $y = w$, perform a plain $\Del$, or $\DS$ with the same $w$ if $w$ is still a neighbor); the resulting function is $\geq f_1$ on $G - u$, so we finish by \zcref{lem:mono}.

\emph{Vertex $x$.} If $x = u$, then $G - x = G - u$ with values $f \geq f_1$; use \zcref{lem:mono}. Otherwise perform the same operation on $G-x$ (a plain $\Del$ if $w = x$); the result dominates $f_1|_{G-u-x}$, and induction plus \zcref{lem:mono} finish.
\end{proof}

\begin{lemma}[postponement]\label{lem:post}
If $v \in V(G)$ satisfies $f(v) > \deg_G(v)$ and $G - v$ is weak-strictly $f$-degenerate, then $G$ is weak-strictly $f$-degenerate.
\end{lemma}
\begin{proof}
Run the sequence for $G-v$ on $G$; each operation lowers $f(v)$ by at most one per removed neighbor, so $v$ keeps value at least $f(v) - \deg_G(v) \geq 1$ and is deleted last.
\end{proof}

\section{From IC-planar graphs to marked plane graphs}
\begin{definition}
A \emph{marked plane graph} is a pair $(H, \F)$ where $H$ is a simple plane graph and $\F$ is a set of inner faces of $H$, each bounded by an induced $4$-cycle of $H$, which are pairwise vertex-disjoint. The faces in $\F$ are called \emph{marked faces}.
For $Q = [a\,c\,b\,d] \in \F$ the \emph{diagonals} of $Q$ are $ab$ and $cd$.
Let $G(H, \F)$ be the simple graph obtained from $H$ by adding both diagonals of every $Q \in \F$.
A vertex $v$ of $H$ is \emph{special} (with respect to $\F$) if it lies on the boundary of some marked face. Since the marked faces are pairwise vertex-disjoint, a special vertex lies on exactly one marked face. We put $\spc(v) = 1$ if $v$ is special and $\spc(v) = 0$ otherwise.
The \emph{boundary} $\pa(H)$ is the set of vertices on the outer face of $H$; the other vertices are \emph{interior}.
\end{definition}

\begin{proposition}\label{prop:diag}
No diagonal of a marked face is an edge of $H$.
\end{proposition}
\begin{proof}
The two ends of a diagonal are opposite vertices of the $4$-cycle bounding the face, and this cycle is induced.
\end{proof}

\begin{lemma}\label{lem:reduction}
For every IC-planar graph $G_0$ there is a marked plane graph $(H, \F)$ with $V(H)=V(G_0)$ and $G_0\subseteq G(H, \F)$.
\end{lemma}
\begin{proof}
Fix an IC-planar drawing of $G_0$ with the minimum number of crossings. Two crossing edges then have four distinct ends (if $ab$ crossed $ac$, exchanging the two initial segments at the crossing removes it without creating new ones, and the drawing stays IC because its crossing set shrinks).
Consider a crossing of $ab$ and $cd$ at a point $p$. The half-curves $ap,pb,cp,pd$ are crossed by nothing else. Delete $ab$ and $cd$ and draw the four curves $a\!-\!c$, $c\!-\!b$, $b\!-\!d$, $d\!-\!a$ as curves running alongside $ap\cup pc$, etc., inside a thin neighborhood of $ap\cup pb\cup cp\cup pd$; these curves cross nothing, and together they bound a face $Q_p=[a\,c\,b\,d]$ whose interior contains only $p$.
If, say, $ac$ was already an edge of $G_0$, that old copy is uncrossed: a crossing on it would share the end $a$ with the crossing at $p$, which IC forbids. So we delete the old copy and keep the new curve, keeping the graph simple.
Do this for every crossing. By IC, the faces $Q_p$ obtained for different crossings are vertex-disjoint, and the procedure for one crossing does not touch the curves used for another.
Each $Q_p$ is bounded by an induced $4$-cycle of $H$: the only edges of $G_0$ joining $a,b$ or $c,d$ are the crossing edges $ab,cd$ themselves ($G_0$ is simple), and these were deleted; the curves added at another crossing join ends of that crossing, which are different from $a,b,c,d$ by IC.
The result is a plane graph $H$ with a set $\F = \{Q_p\}$ of pairwise vertex-disjoint inner $4$-faces bounded by induced $4$-cycles; every edge of $G_0$ is either an edge of $H$ or a diagonal of some $Q_p$ .
\end{proof}

So it suffices to prove that $G(H, \F)$ is weak-strictly $6$-degenerate for every marked plane graph (then use \zcref{lem:sub}).

\section{The inductive statement}\label{sec:inductive}
Let $(H, \F)$ be a marked plane graph and $G = G(H, \F)$. A set $S \subseteq \pa(H)$ is \emph{consecutive} if $S = \emptyset$ or $S$ induces a connected subgraph of the outer face boundary of $H$. For $v \notin S$ put $\deg_S(v) = |N_G(v) \cap S|$ (neighbors in $G$, i.e.\ diagonals count).

\begin{theorem}\label{thm:C}
Let $(H, \F)$ be a marked plane graph, $G = G(H, \F)$, let $S\subseteq\pa(H)$ be a set of at most three consecutive vertices and let $I \subseteq \pa(H) \setminus S$ be independent in $G$. Define
\[
f(v) \coloneqq
\begin{cases}
5+\spc(v) - \deg_S(v), & v \notin \pa(H),\\
4+\spc(v) - \deg_S(v), & v \in \pa(H) \setminus (S \cup I),\\
3-\deg_S(v), & v \in I.
\end{cases}
\]
If $f \geq 1$ (i.e.\ no vertex of $I$ has three $G$-neighbors in $S$), then $G-S$ is weak-strictly $f$-degenerate.
\end{theorem}

\begin{proof}[Proof of \zcref{thm:main}]
Let $G_0$ be an arbitrary IC-planar graph. By \zcref{lem:reduction}, there is a marked plane graph $(H, \F)$ such that $G_0$ is a spanning subgraph of $G(H, \F)$. \zcref{thm:C} with $S = I = \emptyset$ shows that $G(H, \F)$ is weak-strictly $f$-degenerate with $f \leq 6$; \zcref{lem:mono,lem:sub} finish.
\end{proof}

With $\F = \emptyset$, \zcref{thm:C} is \cite[Theorem~2.2]{MR5126524}, restated for weak-strict degeneracy (all values shifted up by one). The only differences are: (i) the diagonals of marked faces are edges of $G$ (and count in $\deg_S$ and in the independence of $I$); (ii) a special vertex not in $I$ gets $+1$; (iii) a vertex of $I$ gets value $3 - \deg_S$ whether it is special or not.

For comparison, we state the planar case.

\begin{theorem}[Bernshteyn, Lee and Smith-Roberge {\cite[Theorem~2.2]{MR5126524}}]\label{thm:BLS}
Let $H$ be a plane graph, let $S \subseteq \pa(H)$ be a set of at most three consecutive vertices, and let $I \subseteq \pa(H) \setminus S$ be independent in $H$. Define
\[
f(v) \coloneqq
\begin{cases}
5 - \deg_S(v), & v \notin \pa(H),\\
4 - \deg_S(v), & v \in \pa(H) \setminus (S \cup I),\\
3 - \deg_S(v), & v \in I.
\end{cases}
\]
If $f \geq 1$ (i.e.\ no vertex of $I$ has three neighbors in $S$), then $H-S$ is weak-strictly $f$-degenerate. Consequently, every planar graph has weak-strict degeneracy at most $5$.
\end{theorem}

\begin{remark}
\zcref{thm:BLS} is exactly the case $\F = \emptyset$ of \zcref{thm:C}, since then $G = H$ and no vertex is special. The three results are related as follows. Yang et al.~\cite{MR4296752} proved that every IC-planar graph is $6$-choosable, using a list-coloring statement of Thomassen type \cite[Theorem~1.2]{MR4296752}. Bernshteyn, Lee and Smith-Roberge proved the weak-strict degeneracy analog of Thomassen's theorem for planar graphs, namely \zcref{thm:BLS}. Our \zcref{thm:C} extends \zcref{thm:BLS} to marked plane graphs, and the resulting \zcref{thm:main} strengthens the theorem of Yang et al.: since $\wsd(G)$ bounds the list chromatic number, $\wsd(G) \leq 6$ implies that $G$ is $6$-choosable, and it also gives the corresponding bounds for DP-coloring and online list coloring. In short,
\[
\text{\zcref{thm:BLS}} \Leftarrow \text{\zcref{thm:C}} \Rightarrow \text{\zcref{thm:main}} \Rightarrow \text{every IC-planar graph is $6$-choosable~\cite{MR4296752}}.
\]
A list-coloring statement such as \cite[Theorem~1.2]{MR4296752} does not by itself give a bound on weak-strict degeneracy, which is why \zcref{thm:main} needs the separate argument of \zcref{sec:proof}.
\end{remark}

\section{Proof of \texorpdfstring{\zcref{thm:C}}{the inductive statement}}\label{sec:proof}

\subsection{Conventions and the deletion rule}
Throughout, $(H, \F)$ is a marked plane graph and $G = G(H, \F)$. We use the following notation.
\begin{itemize}[leftmargin=*,itemsep=1pt]
\item \emph{Adjacency in $G$ and in $H$.} $E(G)=E(H)\cup\{\text{diagonals}\}$. An edge of $H$ is an \emph{$H$-edge}. We write $u \sim v$ if $uv \in E(G)$, and $u \sim_H v$ if $uv\in E(H)$. Then $u$ is a \emph{neighbor} (or \emph{$G$-neighbor}) of $v$, respectively an \emph{$H$-neighbor}; ``$G$-adjacent'' and ``$H$-adjacent'' are used in the same way. $N(v)=N_G(v)$ and $N_H(v)$ are the neighborhoods, and $\deg_G(v)=|N(v)|$, $\deg_H(v)=|N_H(v)|$. For $A\subseteq V(H)$, an \emph{$A$-neighbor} of $v$ is a $G$-neighbor of $v$ in $A$; thus $\deg_S(v)=|N(v)\cap S|$.
\item \emph{Faces.} A face of $H$ is written by listing its boundary vertices in cyclic order, e.g.\ $[a\,b\,c]$ for a triangular face and $[a\,c\,b\,d]$ for a $4$-face. For $[a\,c\,b\,d]\in\F$ the diagonals are $ab$ and $cd$. An \emph{inner face} is a face other than the outer face. For a special vertex $v$, $Q_v$ denotes the marked face containing $v$, and the \emph{partner} $p(v)$ is the vertex opposite to $v$ on $Q_v$. Thus $v$ has exactly one diagonal, namely $vp(v)$. The \emph{rotation} at $v$ is the cyclic order of the $H$-edges at $v$ in the embedding; consecutive edges in it bound a common face.
\item \emph{The outer cycle.} From \zcref{cl:2conn} on, the outer face boundary of $H$ is a cycle $C$. Two vertices of $C$ are \emph{$C$-adjacent} (or \emph{consecutive on $C$}) if they are joined by an edge of $C$, and the two \emph{$C$-neighbors} of $v\in C$ are the vertices $C$-adjacent to it. An \emph{arc} of $C$ is a subpath of $C$.
\item \emph{Other conventions.} $[A]$ is $1$ if the statement $A$ holds and $0$ otherwise (e.g.\ $[y \sim z]$). Inequalities between functions, such as $f' \leq g$, are pointwise.
\end{itemize} A quadruple $(H, \F, S, I)$ satisfying the hypotheses of \zcref{thm:C} (including $f \geq 1$) is an \emph{instance}. Equivalently, $S$ is a set of at most three consecutive boundary vertices, $I \subseteq\pa(H) \setminus S$ is independent in $G$, and $I$ is \emph{admissible}: no vertex of $I$ has three neighbors in $S$. (For $v\notin I$, $f(v) \geq 1$ holds automatically, since $\deg_S(v) \leq 3<4\le\tau(v)$.) An instance is \emph{smaller} than another if it has fewer vertices. Its \emph{type value} is
\[
\tau(v) =
\begin{cases}
5+\spc(v), & v \notin \pa(H),\\
4+\spc(v), & v \in \pa(H) \setminus (S \cup I),\\
3, & v \in I,
\end{cases}
\]
so that $f(v)=\tau(v)-\deg_S(v)$.

\begin{definition}[Deleting a set]\label{def:del}
Let $X \subseteq V(H)$, and assume that no marked face has two opposite vertices in $X$. Let $H_X$ be obtained from $H - X$ as follows. For every $Q = [a\,c\,b\,d] \in \F$ with $Q \cap X = \{a\}$, add the edge $cd$, drawn inside the old face $Q$ close to the path $cbd$. Let $\F_X = \{Q \in \F: Q \cap X = \emptyset\}$.
\end{definition}

The faces of $\F_X$ are still inner $4$-faces of $H_X$, and they are still bounded by induced $4$-cycles, since each new edge $cd$ joins two vertices of a destroyed face, and these lie on no marked face in $\F_X$. So $(H_X,\F_X)$ is a marked plane graph.

\emph{Induced marked faces in auxiliary graphs.} Every other auxiliary graph built below ($H_1^*$ in \zcref{cl:sep}, $H^\circ$ in \zcref{cl:sepcyc}, $H_A,H_B$ in \zcref{lem:PL2}, $H'$ in \zcref{lem:K}) is obtained from a subgraph of $H$ by adding, as an edge, the diagonal of a marked face that is \emph{not} kept in the new marked set. Its ends lie on no remaining marked face, because marked faces are vertex-disjoint, so the remaining marked faces stay induced. The only other added edge is $uw$ in \zcref{cl:tri}, which is checked there. Moreover $G(H_X,\F_X)=G-X$: the only diagonals lost are those with an end in $X$, and the added edges $cd$ are diagonals of $G$. We record the facts about types that we use; all of them are immediate from the construction.

\begin{enumerate}[label=(T\arabic*),leftmargin=*]
\item\label{T1} If $u \in \pa(H) \setminus X$ then $u \in \pa(H_X)$. Indeed, each new edge $cd$ lies inside an old inner face, and it cuts off only the triangle $[c\,b\,d]$.
\item\label{T2} If $u \notin X$ is interior and is $H$-adjacent to a vertex of $X \cap \pa(H)$, then $u \in \pa(H_X)$.
\item\label{T3} If the marked face of $u \notin X$ meets $X$, then $u$ is not special in $H_X$.
\end{enumerate}

\paragraph{The instance after a deletion.} Let $X$ be as in \zcref{def:del} with $X \cap S = \emptyset$, and let $I' \subseteq \pa(H_X) \setminus S$. The quadruple $(H_X, \F_X, S, I')$ has its own type function $\tau'$ and value function $f'$, given by the formulas of \zcref{thm:C} applied to $(H_X, \F_X)$:
\[
\tau'(v) =
\begin{cases}
5+\spc'(v), & v \notin \pa(H_X),\\
4+\spc'(v), & v \in \pa(H_X) \setminus (S \cup I'),\\
3, & v \in I',
\end{cases}
\]
where $\spc'$ refers to $\F_X$, and $f'(v) = \tau'(v) - \deg_S(v)$. Here $\deg_S$ is computed in $G(H_X, \F_X) = G - X$. Since $S \cap X = \emptyset$ and $G - X$ is an induced subgraph of $G$, it equals $\deg_S$ computed in $G$. The quadruple is an instance exactly when $S$ is consecutive in $H_X$ (automatic by \ref{T1}), $I'$ is independent in $G - X$, and no vertex of $I'$ has three neighbors in $S$. The set $I'$ is not determined by $X$: it is chosen in each application, and the lemma below says which choices work.

\begin{lemma}[Deletion rule]\label{lem:R}
Let $(H, \F, S, I)$ be an instance, and let $X \subseteq V(H) \setminus S$ be as in \zcref{def:del}. Consider a sequence of $\Del$/$\DS$ operations on $(G-S, f)$ that removes exactly the vertices of $X$, one at a time. Let $g$ be the resulting function on $G - S - X$, and let $m(u) = f(u) - g(u)$ be the number of times $u$ was lowered. Assume:
\begin{enumerate}[label=(\roman*),leftmargin=*,itemsep=0pt]
\item each $\DS(u, w)$ in the sequence satisfies $f_{\mathrm{current}}(u) > f_{\mathrm{current}}(w)$ for the values current at that moment, and each vertex of $X$ has value at least $1$ at the moment it is removed;
\item $(H_X, \F_X, S, I')$ is an instance for some $I' \subseteq \pa(H_X) \setminus S$;
\item $m(u) \leq \tau(u) - \tau'(u)$ for every $u \in V(H_X) \setminus S$.
\end{enumerate}
If \zcref{thm:C} holds for $(H_X, \F_X, S, I')$, then all operations of the sequence are legal and $G - S$ is weak-strictly $f$-degenerate.
\end{lemma}
\begin{proof}
For $u \in V(H_X) \setminus S$ we have $g(u) = f(u) - m(u) = \tau(u) - \deg_S(u) - m(u)$ and $f'(u) = \tau'(u) - \deg_S(u)$, so (iii) is equivalent to $f' \leq g$.

An operation never raises a value. So, at every moment of the sequence, a vertex $u \in V(H_X) \setminus S$ has value at least its final value $g(u)$, and $g(u) \geq f'(u) \geq 1$ because $(H_X, \F_X, S, I')$ is an instance. A vertex of $X$ has, at every moment before its removal, value at least its value at removal, which is at least $1$ by (i). So all values stay at least $1$ throughout, and with (i) every operation is legal.

By assumption $G - S - X = G(H_X, \F_X) - S$ is weak-strictly $f'$-degenerate, hence weak-strictly $g$-degenerate by \zcref{lem:mono}. Following the sequence by the operations for $(G-S-X, g)$ removes all of $G-S$ by legal operations.
\end{proof}

Since $m \geq 0$, condition (iii) forces $\tau' \leq \tau$ on $V(H_X) \setminus S$. In particular every vertex of $I \setminus X$ must stay in $I'$ (otherwise its type would rise from $3$ to at least $4$), and it may not be lowered at all. So in every application
\[
I' = (I \setminus X) \cup J,
\]
where $J$ is a set of former non-$I$ vertices that we name in each case. A vertex of $J$ has $\tau' = 3$.

When applying \zcref{lem:R} we say that a vertex $u$ \emph{absorbs} its decrements if $m(u) \leq \tau(u)-\tau'(u)$. The typical ways are:
\begin{itemize}[leftmargin=*,itemsep=1pt]
\item an interior vertex that becomes a boundary vertex absorbs $1$;
\item a special vertex whose face is destroyed absorbs $1$;
\item a vertex that enters $I'$ (i.e.\ lies in $J$) absorbs $\tau(u) - 3$.
\end{itemize}
A vertex of $I$ can absorb nothing. In each reduction below we check that no vertex of $I$ is decremented.

\subsection{The minimal counterexample}
Suppose that \zcref{thm:C} fails, and choose a counterexample $(H, \F, S, I)$ that minimizes $|V(H)|$, then maximizes $|E(H)|$, then $|S|$, and then $|I|$. In particular, every instance with fewer vertices satisfies the theorem. Moreover $V(H)\setminus S\neq\emptyset$, since otherwise $G-S$ is the empty graph, which is weak-strictly $f$-degenerate, and $(H, \F, S, I)$ would not be a counterexample.

\begin{claim}\label{cl:conn}
The graph $H$ is connected, and $|S| \geq 2$.
\end{claim}
\begin{proof}
If $H$ is disconnected, we apply the theorem to each component $H_i$, with $S\cap V(H_i)$, $I\cap V(H_i)$ and the marked faces of $H_i$. A marked face has empty interior in $H$, so it is a face of the component that contains it. Passing to a component can only turn interior vertices into boundary vertices, so the values do not increase. Since $S$ is consecutive, it lies in one component, and the resulting sequences of operations can be concatenated.

Now suppose that $|S| \leq 1$, and choose $u \in \pa(H) \setminus S$ as follows: any boundary vertex if $S = \emptyset$, and a neighbor of the vertex of $S$ along the outer boundary walk if $|S| = 1$. Such a vertex exists because $H$ is connected and $V(H) \neq S$. Since $|S \cup \{u\}| \leq 2$, the quadruple $(H, \F, S \cup\{u\}, I \setminus \{u\})$ has values at least $1$, so it is an instance, and its values are exactly those produced by $\Del(G-S,f,u)$. By the maximality of $|S|$ this instance is not a counterexample, and hence neither is $(H, \F, S, I)$.
\end{proof}

\begin{claim}[degrees]\label{cl:deg}
Every interior vertex $v$ has $\deg_G(v) \geq 5 + \spc(v)$, every $v \in \pa(H) \setminus(S \cup I)$ has $\deg_G(v) \geq 4 + \spc(v)$, and every $v \in I$ has $\deg_G(v) \geq 3$. In other words, $\deg_G(v) \geq \tau(v)$ for all $v \notin S$.
\end{claim}
\begin{proof}
Suppose that $\deg_G(v) \leq \tau(v) - 1$, so that $f(v) = \tau(v) - \deg_S(v) > \deg_{G-S}(v)$. Consider the quadruple $(H_{\{v\}}, \F_{\{v\}}, S, I \setminus \{v\})$. By \ref{T1}--\ref{T3} its type values are at most $\tau$, so its values are at most $f$, and by minimality $G - S - v$ is weak-strictly $f$-degenerate. Now \zcref{lem:post} shows that $G-S$ is weak-strictly $f$-degenerate, a contradiction.
\end{proof}

\subsection{Separations and short cycles}
Separations are used only once the outer face boundary of $H$ is a cycle $C$, that is, from \zcref{cl:nochord} on; for a cut vertex, \zcref{cl:2conn} defines $H_1$ and $H_2$ directly. Let $D$ be the closed disc bounded by $C$, so that $H$ is drawn in $D$.

Let $P$ be an induced path of $H$ with two or three vertices, whose ends $x, z$ lie on $C$ and are not consecutive on $C$, and whose middle vertex $y$, if there is one, is interior. As a curve, $P$ joins two points of $C$ through $D$, so it cuts $D$ into two closed discs $D_1, D_2$ with $D_1 \cap D_2 = P$ and $D_1 \cup D_2 = D$. Let $H_i$ be the subgraph of $H$ induced by the vertices lying in $D_i$. We say that $P$ \emph{separates} $H$ into $H_1$ and $H_2$. The following properties hold.
\begin{itemize}[leftmargin=*,itemsep=0pt]
\item $V(H_1) \cup V(H_2) = V(H)$ and $V(H_1) \cap V(H_2) = V(P)$.
\item Every edge of $H$ and every inner face of $H$ lies in $D_1$ or in $D_2$, since none of them can cross the curve $P$. In particular, no edge of $H$ joins $V(H_1) \setminus V(P)$ to $V(H_2) \setminus V(P)$, and the edges of $H_i$ are exactly the edges of $H$ drawn in $D_i$. Indeed, an edge drawn in $D_{3-i}$ with both ends in $D_i$ would have both ends on $P$; when $|V(P)| = 3$ it would be $xz$, which is excluded because $P$ is induced, and when $|V(P)| = 2$ it is $P$ itself.
\item $|V(H_1)| < |V(H)|$ and $|V(H_2)| < |V(H)|$, because each of the two arcs of $C$ between $x$ and $z$ has an inner vertex.
\item The interiors are disjoint: $H_i - V(P)$ lies in $D_i \setminus P$, which is contained in the outer face of $H_{3-i}$. Since $H_i$ is drawn in $D_i$ and $P \subseteq \pa(D_i)$, the path $P$ lies on the outer face boundary of $H_i$, and this boundary is the closed walk formed by $P$ and the arc $C \cap D_i$.
\end{itemize}

Each marked face is an inner face, so it lies in $D_1$ or in $D_2$. Thus $\F = \F_1 \sqcup \F_2$, where $\F_i$ consists of the marked faces inside $H_i$, and no edge of $G$ joins $V(H_1) \setminus V(P)$ to $V(H_2) \setminus V(P)$. The only possible difference between $G[V(H_1)]$ and $G(H_1, \F_1)$ is a diagonal $xz$ of a marked face in $\F_2$ joining the two ends of $P = xyz$. In that case we let $H_1^* = H_1 + xz$, with $xz$ drawn in the outer face of $H_1$ alongside $P$, so that $y$ becomes interior; otherwise we let $H_1^* = H_1$. In both cases $G(H_1^*,\F_1) = G[V(H_1)]$.

\begin{claim}\label{cl:sep}
Let $P$ separate $H$ into $H_1,H_2$, and assume that no vertex of $V(H_1) \setminus S$ has a $G$-neighbor in $S\setminus V(H_1)$, and that $S \cap V(H_1)$ is consecutive in $H_1^*$. For $v\in V(H_2)\setminus(V(P)\cup S)$ let $g(v) = f(v) - |N(v) \cap (V(P) \setminus S)|$. Then either $g(v) < 1$ for some such $v$, or $G[V(H_2)\setminus(V(P)\cup S)]$ is not weak-strictly $g$-degenerate.
\end{claim}
\begin{proof}
The quadruple $(H_1^*,\F_1,S\cap V(H_1),I\cap V(H_1))$ has fewer vertices than $H$, and its type values are at most $\tau$. Indeed, boundary vertices of $H$ stay on the boundary; if $P = xyz$ has three vertices, its middle vertex $y$ is interior in $H$, so its type in $H_1^*$ is at most its type in $H$; and $\spc$ does not increase. By hypothesis, its $\deg_{S\cap V(H_1)}$ equals $\deg_S$. Hence its value function is at most $f$, and on $I \cap V(H_1)$ it equals $f \geq 1$, so the quadruple is indeed an instance. By minimality and \zcref{lem:mono}, $G[V(H_1)]-S$ is weak-strictly $f$-degenerate.

We now run this sequence of operations on $(G-S,f)$. Because $G[V(H_1)]=G(H_1^*,\F_1)$, it acts on $V(H_1)$ exactly as before. Since all $G$-neighbors of a vertex $v \in V(H_2)\setminus (V(P)\cup S)$ outside $V(H_2)\setminus V(P)$ lie in $V(P)$, the sequence lowers $v$ exactly $|N(v)\cap(V(P)\setminus S)|$ times, leaving it with value $g(v)$. If $g \geq 1$, all intermediate values are at least $1$. If moreover $G[V(H_2)\setminus(V(P)\cup S)]$ is weak-strictly $g$-degenerate, then $G-S$ is weak-strictly $f$-degenerate, a contradiction.
\end{proof}

\begin{claim}\label{cl:2conn}
The graph $H$ is $2$-connected. Hence every face boundary is a cycle, and in particular the outer face boundary of $H$ is a cycle $C$.
\end{claim}
\begin{proof}
We have $|V(H)| \geq 3$, since otherwise $V(H) = S$ by \zcref{cl:conn}, which was excluded above. Let $u$ be a cut vertex, let $H_1$ be the subgraph induced by $u$ and one component of $H-u$, and let $H_2$ be the subgraph induced by $u$ and the other components. Then $V(H_1)\cup V(H_2) = V(H)$, $V(H_1) \cap V(H_2) = \{u\}$, no edge of $H$ joins $V(H_1)\setminus\{u\}$ to $V(H_2)\setminus\{u\}$, and $|V(H_i)| < |V(H)|$. Every marked face lies in one $H_i$, since its $4$-cycle lies in one block of $H$. Hence $\F = \F_1 \sqcup \F_2$, $G(H_i,\F_i) = G[V(H_i)]$, and no edge of $G$ joins $V(H_1)\setminus\{u\}$ to $V(H_2)\setminus\{u\}$.

The proof of \zcref{cl:sep} uses only these properties, so \zcref{cl:sep} holds for $P=u$ and these $H_1,H_2$. To apply minimality to an instance on $H_2$, we also need $u$ to lie on the outer face of $H_2$. We use two simple facts. First, a vertex or an edge on the outer face boundary of $H$ stays on the outer face boundary of every subgraph containing it. Second, since $|S| \geq 2$ by \zcref{cl:conn} and $S$ is consecutive, the graph $H[S]$ is connected and contains a path through all vertices of $S$ along the outer face boundary of $H$; this path need not be induced, since for $|S| = 3$ its two ends may be adjacent.

\emph{Case 1: $S \subseteq V(H_i)$ for some $i$.} By relabeling, we may assume that $S\subseteq V(H_1)$. Then $S$ is consecutive in $H_1$, and the hypotheses of \zcref{cl:sep} on $S$ hold trivially. We check that $u$ lies on the outer face of $H_2$. This is clear if $u\in\pa(H)$. Otherwise the outer face boundary of $H$ avoids $u$, so it lies in one component $K_S$ of $H-u$, and $K_S\subseteq H_1$ because it contains $S$. The connected graph $K_S$ is disjoint from $H_2$, so it lies in one face of $H_2$, and this is the outer face, since $K_S$ has vertices on the outer face of $H$. As $u$ is a cut vertex and $H$ is connected, there is an edge from $u$ to $K_S$; it also runs in this face, so $u$ lies on the outer face of $H_2$.

For $v\in V(H_2)\setminus\{u\}$ we have $\deg_S(v)=|N(v)\cap S\cap\{u\}|$, so $g(v)=\tau(v)-[v\sim u]$. The instance $(H_2, \F_2, \{u\}, (I\cap V(H_2)) \setminus \{u\})$ has types at most $\tau$, so its values are at most $g$.

\emph{Case 2: neither $H_1$ nor $H_2$ contains $S$.} Then $u \in S$, since otherwise the connected graph $H[S]$ would lie in one component of $H - u$. As $|S| \leq 3$, we get $S = \{s_1,u,s_2\}$ with $s_1\in V(H_1)\setminus\{u\}$ and $s_2\in V(H_2)\setminus\{u\}$. The vertices $s_1$ and $s_2$ lie in different components of $H-u$, so $s_1s_2\notin E(H)$, and since $H[S]$ is connected, $s_1u$ and $us_2$ are edges on the outer face boundary of $H$. In particular $u\in\pa(H)$, so $u$ lies on the outer face of $H_2$, and $S\cap V(H_1)=\{s_1,u\}$ is consecutive in $H_1$. The only vertex of $S\setminus V(H_1)$ is $s_2$, whose $G$-neighbors lie in $V(H_2)$, so the hypotheses of \zcref{cl:sep} hold, and $g=f$ on $V(H_2)\setminus S$. The instance $(H_2, \F_2, \{u,s_2\}, (I \cap V(H_2)) \setminus S)$ has the consecutive set $\{u,s_2\}$ on the outer face of $H_2$. Since this set has only two vertices, the instance is admissible, and its values are at most $f$.

In both cases, minimality contradicts \zcref{cl:sep}.
\end{proof}

\begin{claim}\label{cl:sepcyc}
No $3$-cycle or $4$-cycle of $H$ has a vertex in its interior.
\end{claim}
\begin{proof}
Let $F$ be such a cycle, let $\mathrm{Int}(F)$ be the set of vertices inside it, let $a \in V(F)$ be arbitrary, and put $S^* = V(F) \setminus \{a\}$. No vertex of $\mathrm{Int}(F)$ is on $\pa(H)$ or in $S$. We first remove the outside of $F$ and then the inside.

For the outside, let $H^\circ$ be obtained from $H - \mathrm{Int}(F)$ by adding, inside $F$, the diagonal of any marked face inside $F$ whose two opposite vertices lie on $F$. There is at most one such face, since two of them would need disjoint opposite pairs of the $4$-cycle $F$, and the first face separates the other pair inside $F$. Let $\F^\circ$ be the set of marked faces outside $F$. Then $G(H^\circ, \F^\circ) = G[V(H) \setminus \mathrm{Int}(F)]$, and the instance $(H^\circ, \F^\circ, S, I)$ is smaller and has values at most $f$. Hence $G - S - \mathrm{Int}(F)$ is weak-strictly $f$-degenerate. Since $N(v) \subseteq \mathrm{Int}(F) \cup V(F)$ for every $v \in \mathrm{Int}(F)$, running that sequence leaves each such $v$ with value
\[
g(v) = f(v) - |N(v) \cap (V(F) \setminus S)| = 5 + \spc(v) - |N(v) \cap V(F)|.
\]

For the inside, let $H^*$ be obtained from the subgraph induced by $\mathrm{Int}(F) \cup V(F)$ by deleting $a$ as in \zcref{def:del}, and let $\F^*$ be the set of marked faces inside $F$ that do not contain $a$. Consider the instance $(H^*, \F^*, S^*, \emptyset)$, whose set $S^*$ is a path on its outer face. For $v \in \mathrm{Int}(F)$ we need $\tau^*(v) - \deg_{S^*}(v) \leq g(v)$, which is equivalent to $\tau^*(v) \leq 5 + \spc(v) - [v\sim a]$. If $va\in E(H)$, then $v$ lies on the outer face of $H^*$ by \ref{T2}, so $\tau^*(v) \leq 4+\spc(v)$. If $v$ is the partner of $a$, then its face is destroyed and $\tau^*(v) \leq 5$. Otherwise $\tau^*(v) \leq 5 + \spc(v)$. So this smaller instance has values at most $g$, and $G - S$ is weak-strictly $f$-degenerate, a contradiction.
\end{proof}

\begin{claim}\label{cl:noK23}
The graph $H$ contains no $K_{2,3}$; that is, no two vertices of $H$ have three common $H$-neighbors. In particular, if $x$ is special, then $x$ and $p(x)$ have no common $H$-neighbor outside $Q_x$.
\end{claim}
\begin{proof}
Suppose that two vertices $u,w$ have three common $H$-neighbors $a,b,c$. The paths $uaw$, $ubw$ and $ucw$ are internally disjoint, and their union has three faces, each bounded by two of the paths. If, say, the outer one is bounded by $uaw$ and $ucw$, then $b$ lies inside the $4$-cycle $uawc$, contradicting \zcref{cl:sepcyc}. For the second statement, let $Q_x=[x\,c_1\,y\,c_2]$ with $y=p(x)$. Then $c_1$ and $c_2$ are common $H$-neighbors of $x$ and $y$, so a third one outside $Q_x$ would give a $K_{2,3}$.
\end{proof}

From now on we write the outer cycle as $C$. An edge of $H$ between two non-consecutive vertices of $C$ is a \emph{chord}.

\begin{claim}\label{cl:nochord}
The cycle $C$ has no chord.
\end{claim}
\begin{proof}
We first consider a chord $uw$ such that, if $|S| = 3$, the middle vertex of $S$ is not in $\{u,w\}$. Since $S$ is consecutive, it lies in one closed side $H_1$, and \zcref{cl:sep} applies with $P=uw$ and $H_1^*=H_1$. Every $v\in V(H_2)\setminus\{u,w\}$ has its $S$-neighbors in $\{u,w\}$, so $g(v)=\tau(v)-|N(v)\cap\{u,w\}|$. This is at least the value function of the smaller instance $(H_2, \F_2, \{u,w\}, (I \cap V(H_2))\setminus\{u,w\})$, which is admissible because its set $S$ has only two vertices, and we obtain a contradiction. Hence we are done if $|S| \leq 2$, and if $S = \{u_1,u_2,u_3\}$ in this order on $C$, then every chord has $u_2$ as an end.

It remains to exclude chords at $u_2$. Let $w_1,\dots,w_k$, where $k \geq 1$, be the other ends of the chords at $u_2$, in order along $C$ from $u_1$ to $u_3$. Put $w_0=u_1$, $w_{k+1}=u_3$ and $W=\{w_1,\dots,w_k\}$. For $0 \leq j \leq k$, let $R_j$ be the subgraph of $H$ induced by the vertices in the closed region between the edges $u_2w_j$ and $u_2w_{j+1}$, and put $S_j=\{w_j,u_2,w_{j+1}\}$, which is consecutive on the outer face of $R_j$. Each marked face lies in one $R_j$, so a vertex of $R_j$ not in $S_j$ has all its $G$-neighbors in $R_j$, and its $G$-neighbors in $S\cup W$ lie in $S_j$.

Let $Z$ be the set of vertices $z \in (I \cap V(R_j)) \setminus S_j$, over all $j$, that are adjacent to all of $S_j$. Such a vertex $z$ is adjacent to $u_2$ by a diagonal, because there is no chord at $u_2$ inside $R_j$; since $u_2$ has at most one diagonal, $|Z| \leq 1$. If $Z=\{z\}$ with $z\in R_j$, then $zw_j$ and $zw_{j+1}$ are $H$-edges, and hence edges of $C$, since every chord ends at $u_2$. So the outer cycle of $R_j$ is the $4$-cycle $u_2w_jzw_{j+1}$, which has no vertex inside by \zcref{cl:sepcyc}. Hence $V(R_j) = \{u_2,w_j,z,w_{j+1}\}$ and $N_G(z) = \{w_j,u_2,w_{j+1}\}$.

We perform the following operations:
\begin{enumerate}[label=(\roman*),leftmargin=*,itemsep=0pt]
\item $\Del(z)$ if $Z=\{z\}$;
\item $\Del(w)$ for every $w\in W\cap I$ (these vertices are pairwise non-adjacent);
\item $\Del(w)$ for the remaining $w\in W$, in any order.
\end{enumerate}

We first check that the vertices of $W$ keep positive values. A vertex $w\in W\cap I$ is not decremented before it is deleted, since its neighbors in $Z\cup W$ that are deleted earlier would have to be in $I$. A vertex $w\in W\setminus I$ is decremented at most $|N(w)\cap(W\cup Z)|$ times. Its neighbors in $S\cup W\cup Z$ are $u_2$, joined to $w$ by the chord, at most two neighbors on $C$, and at most $\spc(w)$ diagonal partners. Therefore,
\[
|N(w) \cap (W \cup Z)| = |N(w) \cap (S \cup W \cup Z)| - |N(w) \cap S| \leq 3 + \spc(w) - |N(w) \cap S| = f(w) - 1,
\]
and the value of $w$ never drops below $1$.

Next consider the remaining vertices. Since $z$ has no neighbor outside $S_j$, after the deletions a vertex $y \in V(R_j) \setminus (S_j \cup Z)$ has value $\tau(y) - |N(y) \cap S_j| - [y \sim z] = \tau(y) - |N(y) \cap S_j|$. For each $j$ with $V(R_j)\ne S_j\cup Z$, consider the instance $(R_j, \F_j, S_j, (I \cap V(R_j)) \setminus S_j)$, where $\F_j$ is the set of marked faces in $R_j$. It is admissible, since an $I$-vertex adjacent to all of $S_j$ would lie in $Z$. Its values are at most the current ones, so in particular the current values are at least $1$. Thus every operation in (i)--(iii) was legal. These instances are smaller, and there are no edges between different graphs $R_j-S_j$, so $G-S$ is weak-strictly $f$-degenerate, a contradiction.
\end{proof}

\begin{claim}\label{cl:notall}
The graph $H$ has an interior vertex.
\end{claim}
\begin{proof}
Otherwise $H = C$, since there are no chords, and $G$ is either $C$ or, if $|C| = 4$ and the inner face is marked, $K_4$. If $G = C$, then every $v \notin S$ has $f(v) \geq 3 - \deg_S(v) = \deg_{G-S}(v) + 1$, and repeated use of \zcref{lem:post} removes $G-S$. If $G = K_4$, then $G - S$ has at most two vertices because $|S| \geq 2$. When there are two, they are adjacent and at most one of them is in $I$; the other one has value $5 - 2 = 3 > 1$ and can be postponed, and the last vertex has value at least $1$.
\end{proof}

\subsection{Faces}

\begin{claim}\label{cl:tri}
Every inner face of $H$ that is not in $\F$ is a triangle.
\end{claim}
\begin{proof}
Let $F = r_1r_2\dots r_k$, with $k \geq 4$, be such a face. If all $r_i$ lie on $C$, then the edges of $F$ are edges of $C$, since there are no chords. In that case $F$ bounds the whole disc of $C$ and $H = C$, contradicting \zcref{cl:notall}.

So we may assume that $r_1 \notin C$. If $r_1r_3 \notin E(H)$, we put $\{u, w\} = \{r_1, r_3\}$. Otherwise the edge $r_1r_3$ runs outside $F$, and one of $r_2$ and $r_4$ lies inside a cycle formed by $r_1r_3$ and part of $F$, thus that vertex is interior. Then $r_2r_4 \notin E(H)$ by planarity, and we put $\{u, w\} = \{r_2, r_4\}$.

Let $z$ be a vertex of $F$ adjacent on $F$ to both $u$ and $w$. We check that $u$ and $w$ are not opposite vertices of a marked face, so that adding $uw$ keeps every marked face induced. Suppose that $Q = [u\,x\,w\,y] \in \F$. If $z \in Q$, then $F$ and $Q$ both contain the path $uzw$, so these are the only faces at $z$. Then $\deg_H(z) = 2$, and $z$ is interior because both faces are inner, so $\deg_G(z) \leq 3$, contradicting \zcref{cl:deg}. Hence $z \notin Q$, and $u$ and $w$ have the three common $H$-neighbors $x$, $y$ and $z$. This is a $K_{2,3}$, which is impossible by \zcref{cl:noK23}.

Adding $uw$ inside $F$ therefore gives a marked plane graph $(H + uw, \F)$. Its types are unchanged, and $I$ is still independent and admissible, because $uw$ has an interior end. Its values are at most $f$, and $H + uw$ has more edges than $H$. By the choice of the counterexample, $G - S + uw$ is weak-strictly $f$-degenerate, and so is its subgraph $G - S$ by \zcref{lem:sub}, a contradiction.
\end{proof}

Consequently, a $4$-cycle of $H$ either bounds a marked face, or has an $H$-chord inside and is split into two triangular faces.

\begin{claim}\label{cl:triang}
We have $|C| \geq 5$ and $|S| = 3$.
\end{claim}
\begin{proof}
By \zcref{cl:notall} there is an interior vertex, so $|C| \geq 5$ by \zcref{cl:sepcyc}.

Suppose that $|S|=2$. Let $u\in C\setminus S$ be a $C$-neighbor of a vertex of $S$, and put $S^+=S\cup\{u\}$. A vertex $v\in C\setminus S^+$ is $H$-adjacent to at most two vertices of $C$, namely its $C$-neighbors, and it has at most one diagonal. If $v$ were adjacent to all three vertices of the path $S^+$, then two of them would be its $C$-neighbors, and these would have to be the two ends of $S^+$, giving $|C|=4$. So $(H,\F,S^+,I\setminus\{u\})$ is an instance. Its values are exactly those produced by $\Del(u)$, which contradicts the maximality of $|S|$.
\end{proof}

\begin{claim}\label{cl:maxI}
Every $v\in C\setminus(S\cup I)$ either has a $G$-neighbor in $I$ or is $G$-adjacent to all three vertices of $S$.
\end{claim}
\begin{proof}
Otherwise $(H,\F,S,I\cup\{v\})$ is an instance whose values are at most $f$, and $|I|$ would not be maximal.
\end{proof}

From now on we write
\[
C=u_1u_2u_3v_1v_2\cdots v_t,\qquad S=\{u_1,u_2,u_3\},\qquad t \geq 2,
\]
and we put $v_{t+1}=u_1$. Since there are no chords, the $H$-neighbors on $C$ of a vertex of $C$ are its two $C$-neighbors, and all its other $H$-neighbors are interior.

\subsection{\texorpdfstring{$2$-chords and pseudo-$2$-chords}{2-chords and pseudo-2-chords}}
A \emph{$2$-chord} is a path $xyz$ of $H$ with $x, z \in C$ not adjacent on $C$ and $y$ interior. Since $C$ has no chord, a $2$-chord is an induced path, so it separates $H$; we let $H_1$ be the side containing $S$ whenever this makes sense.

\begin{lemma}\label{lem:L2}
Let $xyz$ be a $2$-chord with $u_2\notin\{x,z\}$. Then $S$ lies in one side $H_1$, and the other side $H_2$ contains a vertex $a\in I\setminus\{x,z\}$ adjacent in $G$ to $x$, $y$ and $z$. In particular, $x,z\notin I$.
\end{lemma}
\begin{proof}
The set $S$ is a path of $C$ whose middle vertex is not an end of the $2$-chord, so it lies in one closed side $H_1$. \zcref{cl:sep} applies with $P = xyz$, since $S \subseteq V(H_1)$ and the edges of the path $S$ stay on the outer face of $H_1^*$. For $v \in V(H_2) \setminus V(P)$, all $S$-neighbors of $v$ lie in $V(P)$, so $g(v) = \tau(v) - |N(v) \cap V(P)|$.

The instance $(H_2, \F_2, \{x,y,z\}, (I \cap V(H_2)) \setminus \{x,z\})$ has the same types as $H$ on $V(H_2) \setminus V(P)$, so its values equal $g$; a diagonal $xz$ of a marked face in $\F_1$ joins two vertices of the new set $S$ and plays no role. If this instance were admissible, minimality would contradict \zcref{cl:sep}. Hence some vertex of $(I \cap V(H_2)) \setminus \{x,z\}$ is adjacent to all of $x$, $y$ and $z$.
\end{proof}

\begin{lemma}[structure of $2$-chords]\label{lem:TM}
In \zcref{lem:L2}, $V(H_2) = \{x, y, z, a\}$. Moreover, either
\begin{enumerate}[leftmargin=1cm,itemsep=0pt]
\item[(T)] $H_2$ consists of the triangles $xay$ and $ayz$ (so $a$ is not special), or
\item[(M)] $H_2$ is the $4$-cycle $xazy$ bounding a marked face, with diagonals $ay$ and $xz$.
\end{enumerate}
\end{lemma}
\begin{proof}
Let $x = a_0, a_1, \dots, a_k, a_{k+1} = z$ be the arc of $C$ in $H_2$, so that $a = a_i$ for some $1 \leq i \leq k$.

Suppose first that $a$ is joined to $x$ or $z$ by a diagonal, say to $z$. Then $a \neq a_k$, and $ay \in E(H)$, since $a$ has only one diagonal. Since $a$ and $z$ are not adjacent on $C$, the path $ayz$ is a $2$-chord with $u_2 \notin \{a,z\}$, and \zcref{lem:L2} gives $a \notin I$, a contradiction.

Hence $a$ is $C$-adjacent to both $x$ and $z$, and $k=1$. If $ay \in E(H)$, then the triangles $xay$ and $ayz$ are faces by \zcref{cl:sepcyc}, which is case (T).

Otherwise $ay$ is the diagonal of a face $Q_a = [a\,p\,y\,q]$, where $p$ and $q$ are common $H$-neighbors of $a$ and $y$ inside the region bounded by $xazy$. If $x \notin \{p,q\}$, then $a$ and $y$ have the three common $H$-neighbors $p$, $q$ and $x$, forming a $K_{2,3}$, which is impossible by \zcref{cl:noK23}. Hence $x \in \{p,q\}$, say $p=x$, and similarly $q=z$, which is case (M).
\end{proof}

\begin{claim}[diagonals between boundary vertices]\label{cl:diagC}
Let $Q = [a\,c\,b\,d] \in \F$ with $a, b \in C$ and $u_2 \notin \{a,b\}$. Then $a,c,b$ are consecutive on $C$, $c \in I$, $d$ is interior, and $S$ lies on the side of the $2$-chord $adb$ not containing $c$.
\end{claim}
\begin{proof}
By \zcref{prop:diag}, $a$ and $b$ are not adjacent on $C$.

If $c$ and $d$ are both interior, then $acb$ and $adb$ are $2$-chords, and one of them, say $acb$, separates $S$ from $Q$. Its side away from $S$ then contains $a$, $c$, $b$, $d$ and an inner vertex of the other arc of $C$, that is, at least five vertices, contradicting \zcref{lem:TM}.

So we may assume that $c \in C$. Then $ac$ and $cb$ are edges of $C$, and $d$ is interior, since otherwise $|C|=4$. For the $2$-chord $adb$, the side containing $c$ has vertex set $\{a,c,b,d\}$.
\begin{itemize}[leftmargin=*,itemsep=0pt]
\item If $S$ lies on that side, then $S=\{a,c,b\}$, so $\{a, b\} = \{u_1, u_3\}$, and the other side is the side away from $S$. Since $u_2 \notin \{a, b\}$, \zcref{lem:L2} applies to the $2$-chord $adb$, and by \zcref{lem:TM} the other side has exactly one vertex besides $a$, $d$ and $b$. Then $|C| = 4$, contradicting \zcref{cl:triang}. So this case does not occur.
\item Otherwise that side containing $c$ is the side away from $S$. By \zcref{lem:L2,lem:TM}, the vertex $a$ of \zcref{lem:L2} is $c$, so $c\in I$, and (M) holds, since (T) would make $Q$ a pair of triangles. \qedhere
\end{itemize}
\end{proof}

A $2$-chord separates $H$, and \zcref{lem:L2,lem:TM} show that the side away from $S$ is tiny. The next lemma treats a path $awb$ in which $wb$ is a diagonal rather than an edge of $H$. Such a path does not separate $H$, but it still allows a reduction in two phases.

\begin{lemma}[pseudo-$2$-chords]\label{lem:PL2}
There is no interior vertex $w$ together with distinct vertices $a,b \in C\setminus\{u_2\}$, not adjacent on $C$, such that $aw \in E(H)$ and $wb$ is a diagonal of $Q_w = [w\,c_1\,b\,c_2]$, where $a \notin Q_w$ and $c_1,c_2\notin C$.
\end{lemma}
\begin{proof}
Let $A_1$ and $A_2$ be the arcs of $C$ between $a$ and $b$, with $u_2 \in A_1$; then $S\subseteq A_1$, since $u_2\notin\{a,b\}$. The paths $awc_1b$ and $awc_2b$ cut the disc of $C$ into three parts: the face $Q_w$, a part $D_1$ bounded by $A_1$ and one of the paths, and a part $D_2$ bounded by $A_2$ and the other path. By the symmetry between $c_1$ and $c_2$, we may assume that $D_1$ is bounded by $A_1$ and $awc_1b$, and $D_2$ by $A_2$ and $awc_2b$. Let $V_A$ and $V_B$ be the vertex sets of the closed parts $D_1$ and $D_2$. Then $V_A \cap V_B = \{a,w,b\}$, $c_1 \in V_A \setminus V_B$ and $c_2 \in V_B \setminus V_A$. The only edge of $G$ between $V_A \setminus V_B$ and $V_B \setminus V_A$ is the diagonal $c_1c_2$: any other $H$-edge would cross one of the paths, every other marked face lies in $D_1$ or in $D_2$, and $c_1c_2\notin E(H)$ by \zcref{prop:diag}.

In the first phase we remove $V_A\setminus S$. Let $H_A=H[V_A]+wb$, with $wb$ drawn in the outer face next to the path $wc_1b$, so that only $c_1$ becomes enclosed, and let $\F_A$ be the set of marked faces in $D_1$. Then $G(H_A,\F_A)=G[V_A]$. The instance $(H_A,\F_A,S,I\cap V_A)$ is smaller, has $S$ on its outer face, has types at most $\tau$, and has the same $\deg_S$ as $H$. Hence $G[V_A]-S$ is weak-strictly $f$-degenerate.

In the second phase we remove the rest. Running the first phase on $(G-S,f)$ leaves the set $W=V_B\setminus\{a,w,b\}$. Since $S \cap V_B \subseteq\{a,b\}$, the value of $v\in W$ becomes
\[
g(v) = \tau(v) - |N(v) \cap \{a,w,b\}| - [v=c_2].
\]
Let $H_B = H[V_B] + wb$, with $wb$ drawn next to $wc_2b$, so that $c_2$ becomes interior and non-special. Consider the quadruple $(H_B, \F_B, \{a,w,b\}, (I \cap V_B) \setminus \{a,b\})$, whose set $S$ is the path $awb$ on the outer face. For $v \in W \setminus \{c_2\}$ the type is unchanged, and the type of $c_2$ drops from $6$ to $5$, which absorbs the term $[v=c_2]$. So its values are at most $g$, and if it were admissible, minimality would make $G-S$ weak-strictly $f$-degenerate. Hence some $z\in I\cap W$ is adjacent to $a$, $w$ and $b$. The only diagonal at $w$ and the only diagonal at $b$ is $wb$, so $zw$ and $zb$ are $H$-edges. Moreover, $z \in C$, so $z \notin \{c_1,c_2\}$. Thus $w$ and $b$ have the three common $H$-neighbors $c_1$, $c_2$ and $z$, a $K_{2,3}$ in $H$, which is impossible by \zcref{cl:noK23}.
\end{proof}

\begin{lemma}\label{lem:K}
Write $C = u_1u_2u_3v_1 \dots v_t$ as after \zcref{cl:triang}. There is no interior vertex $w$ with $wu_2,wu_3\in E(H)$ such that $wu_1$ is the diagonal of the face $[w\,u_2\,u_1\,c]$. The same holds with $u_1$ and $u_3$ exchanged.
\end{lemma}
\begin{proof}
By \zcref{cl:sepcyc} the triangle $u_2u_3w$ is a face, so $N_H(u_2) = \{u_1, u_3, w\}$, and $c$ is the partner of $u_2$. We first show that $c \notin C$. Since $c \sim_H u_1$, $c \neq u_2$ and there are no chords, $c \in C$ would force $c = v_t$. Then $wv_t \in E(H)$, and $u_3wv_t$ is a $2$-chord, since $u_3$ and $v_t$ are not consecutive on $C$ (as $t \geq 2$) and $u_2$ is not an end. Its side away from $S$ contains the $t+2$ vertices $u_3, w, v_t, v_1, \dots, v_{t-1}$, so $t=2$ by \zcref{lem:TM}. As $w$ already lies on $Q_w$, case (T) holds with $a=v_1$, giving the faces $[u_3\,v_1\,w]$ and $[v_1\,v_2\,w]$. Together with $[u_2\,u_3\,w]$ and $Q_w$ they fill the rotation at $w$, so $N_G(w) = \{u_1,u_2,u_3,v_1,v_2\}$ and $\deg_G(w) = 5 < 6$, contradicting \zcref{cl:deg}. Hence $c$ is interior in $H$.

Let $H'=H-u_2+wu_1$, with $wu_1$ drawn through the old face $Q_w$, and let $S'=\{u_3,w,u_1\}$. The curve $wu_1$ separates $c$ from the new outer face, so $c$ is still interior in $H'$, and the outer cycle of $H'$ is $u_1wu_3v_1\cdots v_t$. We show that $(H',\F\setminus\{Q_w\},S',I)$ is an instance, and then that $\Del(w)$ followed by the sequence for this instance removes $G-S$.

The marked faces in $\F\setminus\{Q_w\}$ stay induced, since the ends of $wu_1$ lie on no other marked face. Also $G(H', \F \setminus\{Q_w\}) = G - u_2$, since $wu_1$ was already a diagonal and the other diagonal $u_2c$ is lost with $u_2$. The set $\{u_3, w, u_1\}$ is a path on the outer cycle of $H'$, the set $I \subseteq\{v_1, \dots, v_t\}$ lies on that cycle, and $I$ is still independent in $G-u_2$. It remains to check that no $z\in I$ is adjacent to all of $S'$. Suppose that $z$ is. Then $z\sim u_1$ and $z\sim w$, and since the only diagonal at $u_1$ and the only diagonal at $w$ is $u_1w$, both $zu_1$ and $zw$ are edges of $H$. Also $z\notin\{u_2,c\}$, since $z\in C\setminus S$ and $c\notin C$. So $u_1$ and $w$ have the three common $H$-neighbors $u_2$, $c$ and $z$ in $H$, a $K_{2,3}$, which is impossible by \zcref{cl:noK23}. Hence the quadruple is an instance.

Next we compare values. After $\Del(w)$, a vertex $y \notin S\cup\{w\}$ has value $\tau(y)-|N(y)\cap S|-[y\sim w]$, while its value in the instance is $\tau'(y)-|N(y)\cap S'|$, where $\tau'\leq\tau$. Since $S\setminus\{u_2\}=S'\setminus\{w\}$, these differ only through $u_2$. As $N_G(u_2)=\{u_1,u_3,w,c\}$, the only affected vertex is $c$, whose type drops from $6$ to $5$; in particular, no vertex of $I$ is affected. So the values after $\Del(w)$ are at least the values $f'$ of the instance.

Minimality applied to this instance shows that $G-S-w$ is weak-strictly $f'$-degenerate, and hence, by \zcref{lem:mono}, weak-strictly degenerate with respect to the values left by $\Del(w)$. These values dominate $f' \geq 1$, so $\Del(w)$ is legal, and $G-S$ is weak-strictly $f$-degenerate. This contradicts the choice of the counterexample.
\end{proof}

\begin{lemma}\label{lem:P}
Let $z\in C\cap I$ and let $x$ be an interior $H$-neighbor of $z$.
\begin{enumerate}[label=(\alph*),leftmargin=*,itemsep=0pt]
\item $x$ has no $H$-neighbor in $I\setminus\{z\}$.
\item If $x$ is special, $z\notin Q_x$, and no $C$-neighbor of $z$ lies on $Q_x$, then $p(x)\notin I$.
\end{enumerate}
\end{lemma}
\begin{proof}
Note that $z\ne u_2$, since $z\notin S$.

(a) Let $b\in I\setminus\{z\}$ with $x\sim_H b$. Then $b\in C\setminus S$, and $b$ is not $C$-adjacent to $z$ because $I$ is independent. So $zxb$ is a $2$-chord whose ends avoid $u_2$, and \zcref{lem:L2} gives $z\notin I$, a contradiction.

(b) Let $Q_x=[x\,c_1\,b\,c_2]$ with $b=p(x)\in I$. Then $b\in C\setminus S$, $b\ne z$, and $b$ is not $C$-adjacent to $z$. Suppose that some $c_i$ lies on $C$. By assumption, $c_i\ne z$ and $c_i$ is not $C$-adjacent to $z$. Also $c_i\ne u_2$, since otherwise $b$, an $H$-neighbor of $u_2$ on $C$, would be in $S$. So $z\,x\,c_i$ is a $2$-chord whose ends avoid $u_2$, and \zcref{lem:L2} gives $z\notin I$, a contradiction. Hence $c_1,c_2\notin C$, and the configuration $z,x,b$ contradicts \zcref{lem:PL2} with $a=z$ and $w=x$.
\end{proof}

\subsection{The final configuration}\label{sec:end}
Recall that $C=u_1u_2u_3v_1\cdots v_t$, $t \geq 2$, $v_{t+1}=u_1$ and $S=\{u_1,u_2,u_3\}$. In this subsection we show that $v_1 \notin I$ (\zcref{cl:E1b}), then that $v_2, v_4 \in I$ (\zcref{cl:E2}), and finally that the inner face at $v_2v_3$ leads to a reducible configuration (\zcref{cl:E3,cl:final}). We use the following facts.
\begin{enumerate}[label=(F\arabic*),leftmargin=*]
\item\label{F1} The inner face of $H$ at an edge $v_iv_{i+1}$ of $C$ is either a marked face, or a triangle $[v_i\,v_{i+1}\,y]$ with $y$ interior, since there are no chords. In the second case we call $y$ the \emph{apex} of $v_iv_{i+1}$. By \zcref{cl:sepcyc} every triangle of $H$ is a face, and $v_iv_{i+1}$ lies on only one inner face, so a common $H$-neighbor of $v_i$ and $v_{i+1}$ must be the apex. In particular, if the inner face at $v_iv_{i+1}$ is in $\F$, then $v_i$ and $v_{i+1}$ have no common $H$-neighbor.
\item\label{F2} A vertex $v_i$ is $G$-adjacent to $u_2$ only through its diagonal. By \zcref{cl:diagC}, if $p(v_i)\in C\setminus\{u_2\}$, then $p(v_i)$ is at distance $2$ from $v_i$ along $C$, the vertex between them is in $I$, and that vertex is not $u_2$. Moreover, neither $v_i$ nor $p(v_i)$ is in $I$. In particular, $p(v_i) \notin I$ for every special $v_i$, since otherwise $p(v_i)$ is interior or $u_2$.
\item\label{F3} A vertex of $C\setminus S$ has at most one diagonal, and its $H$-neighbors on $C$ are its two $C$-neighbors.
\end{enumerate}

By \ref{F2}, $v_1\sim u_1$ only if $t=2$, $[v_1\,v_2\,u_1\,d]\in\F$ and $v_2\in I$.

In each reduction below we apply \zcref{lem:R}, and we only discuss the vertices $u$ with $m(u)>0$ that do not obviously absorb their decrements. Recall that an interior $H$-neighbor of a deleted vertex of $C$ becomes a boundary vertex by \ref{T2}, and that the diagonal partner of a deleted vertex loses its special status by \ref{T3}; either way, such a vertex absorbs one decrement. A vertex of $C$ is decremented by a deleted vertex $v$ only if it is a $C$-neighbor of $v$ or the partner of $v$.

\begin{claim}\label{cl:E1a}
If $v_1\in I$, then $t \geq 3$ and $v_3\in I$. The same holds in the reversed labeling: if $v_t\in I$, then $t \geq 3$ and $v_{t-2}\in I$.
\end{claim}
\begin{proof}
Reversing the labeling of $C$, that is, exchanging $u_1 \leftrightarrow u_3$ and $v_i \leftrightarrow v_{t+1-i}$, exchanges the two statements, so it suffices to prove the first. Suppose that $v_1 \in I$ and that either $t=2$ or $v_3\notin I$. Note that $v_2\notin I$. We apply $\Del(v_1)$, so $X=\{v_1\}$, and we put $I'=(I\setminus\{v_1\})\cup\{v_2\}$.

We first check that every decremented vertex absorbs its decrement. Each neighbor of $v_1$ in $G-S$ is decremented once. The only $C$-neighbor of $v_1$ outside $S$ is $v_2$, and the vertices at distance $2$ from $v_1$ along $C$ are $u_2 \in S$ and $v_3$, where $v_3=u_1\in S$ if $t=2$. So by \ref{F2}, either $p(v_1) \notin C$, or $p(v_1)=u_2$, or $p(v_1)=v_3$ and the vertex $v_2$ between them is in $I$. The last case is impossible, since $v_1 \in I$ and $I$ is independent. If $p(v_1)=u_2$, then $p(v_1) \in S$ and there is nothing to check; otherwise $p(v_1)$ is interior and absorbs its decrement by \ref{T3}. The vertex $v_2$ enters $I'$, so it absorbs $1+\spc(v_2) \geq 1$. No vertex of $I$ is a neighbor of $v_1$. Since there are no chords, the remaining neighbors of $v_1$ are its $H$-neighbors off $C$, which are interior in $H$. Each of them is $H$-adjacent to $v_1\in X\cap\pa(H)$, so it becomes a boundary vertex of $H_X$ by \ref{T2} and absorbs its decrement.

It remains to check that $I'$ is admissible. The neighbors of $v_2$ are $v_1$, the vertex $v_3$ (which is in $S$ if $t=2$ and not in $I$ otherwise), interior vertices, and $p(v_2)$. By \ref{F2}, $p(v_2)\in C\setminus S$ would mean $p(v_2)=v_4$ with $v_3\in I$, which is excluded. So $v_2$ has no neighbor in $I\setminus\{v_1\}$. Finally, $v_2$ is adjacent to $u_2$ and to $u_3$ only through its diagonal, so $\deg_S(v_2) \leq 2$.

By minimality and \zcref{lem:R}, $G-S$ is weak-strictly $f$-degenerate, a contradiction.
\end{proof}

\paragraph{Orientation.} If $u_2$ is special, then $p(u_2)$ is not in $\{u_1,u_3\}$ by \zcref{prop:diag}, so $p(u_2)=v_i$ for at most one $i$. Reversing the labeling of $C$ if necessary, we assume that
\[
p(u_2)\ne v_1.
\]
Everything before this point is symmetric under the reversal, and from now on the labeling is fixed. By \ref{F2}, the orientation gives $v_1 \not \sim u_2$. It is used directly only in the proof of \zcref{cl:E1b}; \zcref{cl:E2,cl:E3,cl:final} use it only through \zcref{cl:E1b}.

\begin{claim}\label{cl:E1b}
The vertex $v_1$ is not in $I$.
\end{claim}
\begin{proof}
Suppose that $v_1\in I$. By \zcref{cl:E1a}, $t \geq 3$ and $v_3\in I$, so $v_2\notin I$. Since $t \geq 3$, the orientation gives $\deg_S(v_1)=1$, so $f(v_1)=2$. Moreover, $\deg_S(v_2) \leq 1$, with equality only through the diagonal of $v_2$, in which case $\spc(v_2)=1$. Hence $f(v_2) \geq 4 > 3 \geq f(v_3)$.

We apply $\DS(v_2, v_3)$ and then $\Del(v_1)$, so $X=\{v_1, v_2\}$. The first operation lowers $v_1$ to $1$, which is legal, and $v_1$ is deleted next. By \zcref{prop:diag}, no face has two opposite vertices in $X$.

We first count the decrements on $C \setminus S$. The vertex $v_3$ is saved, and $v_3\not\sim v_1$ since both are in $I$, so $m(v_3)=0$. The other vertices of $C \setminus S$ that may be decremented are the partners of $v_1$ and $v_2$. By \ref{F2}, $p(v_1)\notin C\setminus S$, since that would need $v_1 \sim v_3$. Also by \ref{F2}, $p(v_2)\in C\setminus S$ means $p(v_2)=v_4$, since $v_3\in I$; in that case $v_4 \notin I$, and it absorbs by \ref{T3}. No vertex of $I$ is decremented, since the $I$-neighbors of $v_2$ are $v_1$ and $v_3$.

By \ref{F3}, this covers all vertices of $C\setminus S$, so it remains to consider interior vertices $y$. If $y$ is adjacent to only one of $v_1$ and $v_2$, it absorbs its decrement by \ref{T2} or \ref{T3}. Suppose that $y$ is adjacent to both. If $Q_y$ meets $X$, then $y$ loses its special status and, being $H$-adjacent to at least one vertex of $X$, becomes a boundary vertex, so it absorbs $2$. Otherwise either $y$ is not special or $Q_y \cap \{v_1,v_2\} = \emptyset$, so $y$ is $H$-adjacent to both $v_1$ and $v_2$, and by \ref{F1} it is the apex $w$ of $v_1v_2$. In that case we put $I' = (I \setminus\{v_1\}) \cup \{w\}$, so that $w$ absorbs $\tau(w) - 3 \geq 2$; in all other cases we put $I' = I \setminus \{v_1\}$.

It remains to check that $I'$ is admissible, which concerns only $w$. By \ref{T2}, $w$ is a boundary vertex of $H_X$. By \zcref{lem:P}(a) with $z=v_1$, it has no $H$-neighbor in $I\setminus\{v_1\}$. Next let $b = p(w)$, and suppose that $b \in I \setminus \{v_1\}$. If $u_3\in Q_w$, then $b$ is an $H$-neighbor of $u_3$ different from $v_1$ and $u_2$, so $u_3b$ is a chord of $C$, a contradiction. Otherwise no $C$-neighbor of $v_1$ lies on $Q_w$, and \zcref{lem:P}(b) with $z=v_1$ gives $b\notin I$.

Finally, suppose that $w$ is adjacent to $u_1$, $u_2$ and $u_3$. If $w$ is $H$-adjacent to $u_1$, then $u_1wv_2$ is a $2$-chord avoiding $u_2$ (recall that $t \geq 3$), and its side not containing $S$ contains $u_1,w,v_2,v_3,\dots,v_t$. By \zcref{lem:TM} this forces $t=3$, and either (M) holds, in which case $v_2\in Q_w$, which is excluded, or (T) holds with the faces $[u_1\,v_3\,w]$ and $[v_3\,w\,v_2]$. In the latter case the faces at $v_2$ are $[v_1\,v_2\,w]$, $[v_2\,v_3\,w]$ and the outer face, so $\deg_G(v_2) = 3 < 4$, contradicting \zcref{cl:deg}. Hence $wu_1$ is a diagonal of $Q_w$, and $wu_2,wu_3\in E(H)$. Since $u_2$ is a common $H$-neighbor of $w$ and $u_1=p(w)$, \zcref{cl:noK23} gives $u_2\in Q_w$, so $Q_w=[w\,u_2\,u_1\,c]$, which contradicts \zcref{lem:K}.

So $\deg_S(w) \leq 2$, and $I'$ is independent and admissible. Now \zcref{lem:R} and minimality give a contradiction.
\end{proof}

\begin{claim}\label{cl:E2}
We have $v_2\in I$, $v_3\notin I$, $t \geq 4$ and $v_4\in I$.
\end{claim}
\begin{proof}
By \zcref{cl:E1b,cl:maxI}, $v_1$ has a $G$-neighbor in $I$ or is adjacent to all of $S$. The latter is impossible, since $v_1$ is not $C$-adjacent to $u_1$ or $u_2$ and it has at most one diagonal. The neighbors of $v_1$ on $C \setminus S$ are $v_2$ and possibly $p(v_1)$, and $p(v_1) \notin I$ by \ref{F2}. Hence $v_2 \in I$, and then $v_3 \notin I$. If $t=2$, then $v_2=v_t\in I$, and the reversed form of \zcref{cl:E1a} gives $t \geq 3$. So $t \geq 3$.

Now suppose that $v_4 \notin I$; this includes the case $t=3$, in which $v_4=u_1$. We apply $\Del(v_2)$, so $X=\{v_2\}$. The neighbors of $v_2$ in $G - S - X$ are $v_1$, $v_3$, interior vertices and $p(v_2)$. The vertex $p(v_2)$ absorbs its decrement by \ref{T3}, and it is not in $I$, because by \ref{F2}, $p(v_2)\in C\setminus S$ would need $v_1\in I$ or $v_3\in I$. In particular, no vertex of $I$ is decremented. We choose $I'$ according to whether $v_2$ lies on a marked face $[v_1\,v_2\,v_3\,d]$.

If $[v_1\,v_2\,v_3\,d]\in\F$, we put $I'=I\setminus\{v_2\}$; then $v_1$ and $v_3$ absorb their decrements by \ref{T3}.

Otherwise $v_1 \not\sim v_3$ by \ref{F2}, and we put $I' = (I \setminus \{v_2\}) \cup \{v_1,v_3\}$, so that $v_1$ and $v_3$ absorb their decrements by entering $I'$. They are independent of $I \setminus \{v_2\}$. The $C$-neighbors of $v_1$ and $v_3$ other than $v_2$ are $u_3$ and $v_4$, both outside $I$. By \ref{F2}, $\{p(v_1), p(v_3)\} \cap I = \emptyset$. Finally, $\deg_S(v_1) \leq 2$, since $v_1$ is adjacent to $u_3$ and has at most one diagonal, and $\deg_S(v_3) \leq 2$, since $v_3$ is adjacent to $u_1$ only if $t=3$ and has at most one diagonal.

In both cases \zcref{lem:R} and minimality give a contradiction. Hence $v_4\in I$; in particular, $v_4\ne u_1$, so $t \geq 4$.
\end{proof}

From now on, $v_1,v_3\notin I$, $v_2,v_4\in I$, and $t \geq 4$. Hence $\deg_S(v_2) \leq 1$ and $\deg_S(v_3) \leq 1$, since both vertices are adjacent to $S$ only through their diagonals, and so
\[
f(v_2)=3-\deg_S(v_2) \geq 2,\qquad f(v_3)=4+\spc(v_3)-\deg_S(v_3)\geq 4 > 3 \geq f(v_4).
\]

\begin{claim}\label{cl:E3}
The inner face at $v_2v_3$ is a triangular face $[v_2\,v_3\,x]$ such that $x \sim v_1$, and either $x$ is not special or $Q_x \cap \{v_2,v_3\} = \emptyset$.
\end{claim}
\begin{proof}
Suppose that the claim fails. We apply $\DS(v_3,v_4)$ and then $\Del(v_2)$, so $X=\{v_2,v_3\}$. The first operation lowers $v_2$ to $f(v_2)-1 \geq 1$, which is legal. We will put $v_1$ into $I'$ in every case, and possibly also the apex $x$ of $v_2v_3$.

We first count the decrements. The vertex $v_4$ is saved, and $v_4 \not\sim v_2$ since $\{v_2, v_4\} \subseteq I$, so $m(v_4)=0$. The vertex $v_1$ is decremented by $v_2$, and also by $v_3$ if $[v_1\,v_2\,v_3\,d]\in\F$; in the latter case $\spc(v_1)=1$. Since $v_1$ enters $I'$, it absorbs $1+\spc(v_1) \geq m(v_1)$. The partners of $v_2$ and $v_3$ absorb their decrements by \ref{T3}, and they are not in $I$ by \ref{F2}, so no vertex of $I$ is decremented. As in the proof of \zcref{cl:E1b}, a common neighbor $y\notin\{v_1,v_4\}$ of $v_2$ and $v_3$ is interior by \ref{F3}, and it absorbs $2$ if $Q_y$ meets $X$; otherwise $y=x$ is the apex of $v_2v_3$.

The vertex $v_1$ has no neighbor in $I\setminus\{v_2\}$, since its $C$-neighbors are $u_3$ and $v_2$, and by \ref{F2}, neither $v_1$ nor $p(v_1)$ is in $I$. Hence $(I\setminus \{v_2\}) \cup \{v_1\}$ is independent. Moreover, $\deg_S(v_1) \leq 2$, since $v_1$ is adjacent to $u_3$ and has at most one diagonal.

We now define $I'$ and consider the quadruple $(H_{X}, \F_{X}, S, I')$. If the inner face at $v_2v_3$ is in $\F$, or it is a triangular face $[v_2\, v_3\, x]$ and $Q_x$ meets $X$, we put $I'=(I\setminus\{v_2\})\cup\{v_1\}$. Otherwise the inner face at $v_2v_3$ is a triangular face $[v_2\, v_3\, x]$ such that either $x$ is not special or $Q_x \cap \{v_2,v_3\} = \emptyset$, and since the claim fails, $x \not\sim v_1$. In this case we put $I' = (I \setminus \{v_2\}) \cup \{v_1,x\}$, so that $x$ absorbs $\tau(x)-3 \geq 2$. In the second case we must check that $x$ can join $I'$. By \ref{T2}, $x$ lies on the boundary of $H_X$, and by \zcref{lem:P}(a) with $z=v_2$, it has no $H$-neighbor in $I\setminus\{v_2\}$. Since $x \not\sim v_1$, we have $v_1 \notin Q_x$, so no $C$-neighbor of $v_2$ lies on $Q_x$, and \zcref{lem:P}(b) with $z=v_2$ gives $p(x)\notin I$. To see that $\deg_S(x) \leq 2$, suppose first that $x\sim_H u_3$. Then the $2$-chord $u_3xv_3$ avoids $u_2$, and its side away from $S$ contains $u_3,x,v_3,v_1,v_2$, contradicting \zcref{lem:TM}. Similarly, if $x\sim_H u_1$, then the $2$-chord $u_1xv_2$ avoids $u_2$ (note that $u_1$ and $v_2$ are not $C$-adjacent, as $t \geq 4$), and its side away from $S$ contains $u_1,x,v_2,v_3,v_4$, again a contradiction. So $x$ is adjacent to at most one of $u_1$ and $u_3$, through its diagonal, and $\deg_S(x)\leq 2$. Together with $x\not\sim v_1$, this shows that $I'$ is independent and admissible.

Thus in both cases $(H_{X}, \F_{X}, S, I')$ is an instance, $I'$ is independent and admissible, and every vertex in $N_{G}(X) \setminus S$ absorbs its decrements. Now \zcref{lem:R} and minimality give a contradiction.
\end{proof}

\begin{claim}\label{cl:final}
The configuration of \zcref{cl:E3} does not occur.
\end{claim}
\begin{proof}
Let $x$ be as in \zcref{cl:E3}. If $x \sim v_1$ by a diagonal, then $p(x)=v_1$, and $v_2\notin Q_x$ is a common $H$-neighbor of $x$ and $v_1$, which contradicts \zcref{cl:noK23}.

So $xv_1\in E(H)$. Then $v_1v_2x$ is a triangle, and by \ref{F1} it is the inner face at $v_1v_2$. The faces at $v_2$ are therefore $[v_1\,v_2\,x]$, $[v_2\,v_3\,x]$ and the outer face. So $N_G(v_2)=\{v_1,x,v_3\}$, the vertex $v_2$ is not special, $\deg_S(v_2)=0$ and $f(v_2)=3$. In particular, $[v_1\,v_2\,v_3\,d]\notin\F$, so $v_1\not\sim v_3$. Moreover, $f(v_1)\geq 3$, since the only $C$-neighbor of $v_1$ in $S$ is $u_3$, so that $\deg_S(v_1) \leq 1+\spc(v_1)$.

We apply $\DS(v_3,v_4)$ and then $\DS(v_1,v_2)$, and we let $g$ be the resulting value function on $G - S - \{v_1, v_3\}$. The first operation is legal as before; it lowers $v_2$ to $2$ and does not touch $v_1$. The second one is defined because $f(v_1) \geq 3>2$. Now $v_2$ has value $2$, and its only remaining neighbor is $x$, so its value is larger than its degree. Let $X = \{v_1,v_2,v_3\}$. By \zcref{lem:post}, it suffices to show that $G-S-X$ is weak-strictly $g$-degenerate; note that $v_2$, although it lies in $X$, does not lower $x$. As in the proof of \zcref{lem:R}, it suffices to find an instance $(H_X,\F_X,S,I')$ with values $f'\le g$, that is, with $m(u)\le\tau(u)-\tau'(u)$, where $m(u)$ counts only the decrements caused by $v_1$ and $v_3$. No face has two opposite vertices in $X$, since $v_1$ and $v_3$ would have to be opposite on a face containing $v_2$. We put $I'=(I\setminus\{v_2\})\cup\{x\}$.

We first count the decrements. As before, $m(v_4)=0$. The vertex $x$ is $H$-adjacent to $v_1$ and $v_3$, so $m(x)=2$, and $x$ absorbs $\tau(x)-3 \geq 2$. The partners of $v_1$ and $v_3$ absorb their decrements by \ref{T3}, and they are not in $I$ by \ref{F2}. A vertex $y\ne x$ adjacent to only one of $v_1$ and $v_3$ absorbs its single decrement. Now let $y\ne x$ be adjacent to both; then $y\notin C$ by \ref{F3}. If $y\sim_H v_1$ and $y\sim_H v_3$, then $v_1yv_3$ is a $2$-chord avoiding $u_2$, and \zcref{lem:TM} applies with $a=v_2$. Here (M) is impossible since $v_2$ is not special, and (T) gives $y\sim_H v_2$, so $y=x$. Hence one of the two adjacencies is a diagonal, $Q_y$ meets $X$, and $y$ absorbs $2$ by \ref{T2} and \ref{T3}. No vertex of $I$ is decremented, since the $I$-neighbors of $v_1$ and $v_3$ are $v_2$ and $v_4$, and $v_4$ is saved.

It remains to check that $I'$ is admissible. By \ref{T2}, $x$ is a boundary vertex of $H_X$, and by \zcref{lem:P}(a), it has no $H$-neighbor in $I\setminus\{v_2\}$. For $p(x)$, recall from \zcref{cl:E3} that $v_3\notin Q_x$. If $v_1\in Q_x$, then $p(x)$ is an $H$-neighbor of $v_1$ other than $v_2$, and hence not in $I$; otherwise \zcref{lem:P}(b) with $z=v_2$ applies. Finally, $\deg_S(x) \leq 2$ by the argument in the proof of \zcref{cl:E3}.

By minimality, $G-S-X$ is weak-strictly $f'$-degenerate, and hence weak-strictly $g$-degenerate. Postponing $v_2$, we conclude that $G-S$ is weak-strictly $f$-degenerate, a contradiction.
\end{proof}

\zcref{cl:E2,cl:E3,cl:final} are contradictory, so there is no counterexample. This proves \zcref{thm:C}.

\section*{Acknowledgments}
This project was initiated about three years ago, when Zhen He wrote a draft containing some gaps. Recently, we used Claude Opus 5.5 to fill these gaps, in particular, to establish the form of \zcref{thm:C}. We have read all the content generated by Claude and have revised and edited it, and we take full responsibility for the correctness and content of this manuscript.
\section*{Conflict of Interest Statement}
The authors declare no conflicts of interest.

\section*{Data Availability Statement}
Data sharing is not applicable to this article as no new data were created or analyzed in this study.

\printbibliography
\end{document}

%% file: Preamble.tex
\usepackage{mathtools}
\usepackage{iftex}

\ifPDFTeX
  \usepackage[T1]{fontenc}   % Better font encoding (supports accented chars)
  \usepackage[utf8]{inputenc} % UTF-8 input support
  \usepackage{lmodern}       % Latin Modern font family
  \usepackage{amssymb}       % \mathbb etc. (unicode-math provides these otherwise)
\else
  \usepackage{fontspec}      % Load OpenType fonts (Unicode)
  \usepackage[warnings-off={mathtools-overbracket,mathtools-colon}]{unicode-math}
\fi

\ifXeTeX
  \usepackage[
    protrusion=true,
    expansion=false,
    final,
    factor=1100,
  ]{microtype}
\else
  \usepackage[
    activate={true,nocompatibility},
    final,
    tracking=true,
    factor=1100,
  ]{microtype}
\fi

\usepackage[a4paper,margin=25mm]{geometry} % Uniform margins
\usepackage{setspace}
\usepackage{amsthm}
\allowdisplaybreaks[4] % Allow page breaks in align environments

\usepackage[inline]{enumitem} % Flexible list control

\usepackage{graphicx} % Include images
\usepackage{xcolor}   % Color support
\usepackage{tikz}     % Drawing diagrams
\usetikzlibrary{
    arrows.meta,
    calc,
    intersections,
    positioning,
    shapes.geometric,
    decorations.pathreplacing,
    decorations.markings
}

\usepackage{caption}
\usepackage[labelformat=simple]{subcaption}

\usepackage[
  backend=biber,            % Use biber (required for Unicode)
  citestyle=numeric-comp,       % Compressed numeric citations [1–3]
  bibstyle=WT-numeric,
  sorting=nty,              % Name–Title–Year sorting
  maxbibnames=99,
  maxcitenames=99,
  isbn=false,
  url=false
]{biblatex}

\AtBeginBibliography{%
  \ifXeTeX
    \microtypesetup{protrusion=true}%
  \else
    \microtypesetup{protrusion=true, expansion=true}%
  \fi
}

\newif\ifshowdoi
\showdoitrue   % <-- default: show DOI

\usepackage[
  colorlinks=true,
  linkcolor=blue,
  citecolor=black,
  urlcolor=blue
]{hyperref}

\usepackage{zref-clever}
\zcsetup{cap=true, abbrev=true}
\zcRefTypeSetup{claim}{Name-sg=Claim,Name-pl=Claims}
\zcRefTypeSetup{case}{Name-sg=Case,Name-pl=Cases}
\zcRefTypeSetup{conjecture}{Name-sg=Conjecture,Name-pl=Conjectures}
\zcRefTypeSetup{problem}{Name-sg=Problem,Name-pl=Problems}
\zcRefTypeSetup{observation}{Name-sg=Observation,Name-pl=Observations}
\makeatletter
\renewenvironment{proof}[1][\proofname]{\par
  \pushQED{\qed}%
  \normalfont
  \topsep6\p@\@plus6\p@\relax
  \trivlist
  \item[\hskip\labelsep
        \bfseries
        #1\@addpunct{.}]\ignorespaces
}{%
  \popQED\endtrivlist\@endpefalse
}
\makeatother

\theoremstyle{definition}
\newtheorem{theorem}{Theorem}[section]
\newtheorem{lemma}[theorem]{Lemma}

\newtheorem{claim}{Claim}

\newtheorem{definition}{Definition}
\newtheorem{remark}{Remark}
\newtheorem{proposition}{Proposition}

\renewcommand{\theHtheorem}{\thesection.\arabic{theorem}}

\makeatletter
\newcommand{\sharetheH}[1]{%
  \expandafter\let\csname theH#1\endcsname\theHtheorem
}
\makeatother

\sharetheH{lemma}
\sharetheH{corollary}

\AddToHook{env/lemma/begin}{\zcsetup{countertype={theorem=lemma}}}
\AddToHook{env/corollary/begin}{\zcsetup{countertype={theorem=corollary}}}

\tikzset{
  mid arrow/.style={
    postaction={
      decorate,
      decoration={
        markings,
        mark=at position .55 with {\arrow{stealth}}
      }
    }
  }
}

\tikzset{
  hollow node/.style={
    circle, draw, fill=white, inner sep=1pt
  },
  solid node/.style={
    circle, draw, fill=black, inner sep=1pt
  },
  rect solid node/.style={
    rectangle, draw, fill=green, inner sep=2pt
  }
}